\documentclass{amsart}

\usepackage{fancyhdr}
\usepackage{lastpage}
\usepackage{stmaryrd,yhmath}

\newcommand{\bdis}{\begin{displaymath}}
\newcommand{\edis}{\end{displaymath}}
\newcommand{\be}{\begin{equation}}
\newcommand{\ee}{\end{equation}}

\newcommand{\mcal}{\mathcal}

\theoremstyle{definition}

\newtheorem{cor}[]{Corollary}

\theoremstyle{remark}
\newtheorem{remark}[]{Remark}

\newtheorem*{mydef1}{{\bf Theorem}}

\numberwithin{equation}{section}

\begin{document}

\title{Jacob's ladders and the first asymptotic formula for the expression of the fifth order
$Z[\varphi(t)/2+\rho_1]Z[\varphi(t)/2+\rho_2]Z[\varphi(t)/2+\rho_3]\hat{Z}^2(t)$ for the collection of disconnected sets}

\author{Jan Moser}

\address{Department of Mathematical Analysis and Numerical Mathematics, Comenius University, Mlynska Dolina M105, 842 48 Bratislava, SLOVAKIA}

\email{jan.mozer@fmph.uniba.sk}


\begin{abstract}
It is shown in this paper that there is a fine correlation of the fifth order between the values
$Z[\varphi(t)/2+\rho_1]Z[\varphi(t)/2+\rho_2]Z[\varphi(t)/2+\rho_3]$ and $\hat{Z}^2(t)$ which correspond to two collections of disconnected sets.
This new asymptotic formula cannot be obtained within known theories of Balasubramanian, Heath-Brown and Ivic.
\end{abstract}

\maketitle

\section{The main result}

\subsection{}

Let (see \cite{4}, p. 24; $G_3\to G_5,\ G_4\to G_6$)
\begin{eqnarray}
& &
G_5(x_1,x_2)=G_5(x_1,x_2,T,U)= \nonumber \\
& &
\bigcup_{T\leq k_{2\nu}\leq T+U}\left\{ t:\ k_{2\nu}(x_1)\leq t\leq k_{2\nu}(x_2)\right\}, \nonumber \\
& &
x_1,x_2\in [-\pi/2,\pi/2],\ x_1<x_2 , \nonumber \\
& &
G_6(y_1,y_2)=G_6(y_1,y_2,T,U)= \\
& &
\bigcup_{T\leq k_{2\nu+1}\leq T+U}\left\{ t:\ k_{2\nu+1}(y_1)\leq t\leq k_{2\nu+1}(y_2)\right\}, \nonumber  \\
& &
y_1,y_2\in [-\pi/2,\pi/2],\ y_1<y_2 , \nonumber
\end{eqnarray}
and the sequence $\{ k_\nu(\tau)\}$ is defined by the equation (see \cite{5}, p. 136)
\bdis
\vartheta_1[k_\nu(\tau)]=\frac{1}{3}(\pi\nu+\tau),\ \nu=1,2,\dots  ,\ \tau\in [-\pi,\pi] ,
\edis
where
\bdis
\vartheta_1(t)=\frac{t}{2}\ln\frac{t}{2\pi}-\frac{t}{2}-\frac{\pi}{8},\ k_\nu=k_\nu(0) .
\edis

\begin{remark}

Collections of disconnected sets $G_1(x),G_2(y)$ and $G_3(x), G_4(y)$ are connected with the sequences $\{ t_\nu(\tau)\},\ \{ g_\nu(\tau)\}$ which
have been defined in \cite{2} and \cite{3},\cite{4}.

\end{remark}

Let the values
\begin{eqnarray*}
& &
T,\ k_{2\nu}(x_1),\ k_{2\nu}(x_2),\ k_{2\nu+1}(y_1),\ k_{2\nu+1}(y_2),\ T+U \\
& &
\mathring{T},\ \mathring{k}_{2\nu}(x_1),\ \mathring{k}_{2\nu}(x_2),\ \mathring{k}_{2\nu+1}(y_1),\ \mathring{k}_{2\nu+1}(y_2),\ \widering{T+U}
\end{eqnarray*}
correspond each other by following equations

\be \label{e1.2}
T=\frac{1}{2}\varphi(\mathring{T}),\ k_{2\nu}(x_1)=\frac{1}{2}\varphi[\mathring{k}_{2\nu}(x_1)],\dots , T+U=\frac{1}{2}\varphi[\widering{T+U}]
\ee
($\varphi(t),\ t\geq T_0[\varphi]$ is an increasing function). Thus the mapping of (1.1)
\be \label{e1.3}
G_5(x_1,x_2)\rightarrow\mathring{G}_5(x_1,x_2),\ G_6(y_1,y_2)\rightarrow\mathring{G}_6(y_1,y_2)
\ee
is defined.

\subsection{}

The following theorem holds true

\begin{mydef1}

\begin{eqnarray}
& &
\int_{\mathring{G}_5(x_1,x_2)}Z\left[\frac{\varphi(t)}{2}+\rho_1\right]Z\left[\frac{\varphi(t)}{2}+\rho_2\right]
Z\left[\frac{\varphi(t)}{2}+\rho_3\right]\hat{Z}^2(t){\rm d}t=\nonumber \\
& &
\frac{4}{\pi}U\sin\frac{x_2-x_1}{2}\cos\left\{\frac{x_1+x_2}{2}+(\rho_1+\rho_2+\rho_3)\ln P\right\}+\mcal{O}(T^{13/16+\epsilon}), \nonumber \\
& & \ \\
& &
\int_{\mathring{G}_6(y_1,y_2)}Z\left[\frac{\varphi(t)}{2}+\rho_1\right]Z\left[\frac{\varphi(t)}{2}+\rho_2\right]
Z\left[\frac{\varphi(t)}{2}+\rho_3\right]\hat{Z}^2(t){\rm d}t=\nonumber \\
& &
-\frac{4}{\pi}U\sin\frac{y_2-y_1}{2}\cos\left\{\frac{y_1+y_2}{2}+(\rho_1+\rho_2+\rho_3)\ln P\right\}+\mcal{O}(T^{13/16+\epsilon}), \nonumber
\end{eqnarray}
where
\be \label{e1.5}
t-\frac{\varphi(t)}{2}\sim (1-c)\pi(t),\ t\to\infty ,
\ee
\be \label{e1.6}
T^{13/16+2\epsilon}\leq U\leq T^{7/8+\epsilon/2},\ \rho_1,\rho_2,\rho_3=\mcal{O}(T^{1/48-\epsilon}),\ P=\sqrt{\frac{T}{2\pi}} ,
\ee
$c$ is the Euler's constant and $\pi(t)$ is the prime-counting function.

\end{mydef1}

Since (see (1.2), (1.5))
\bdis
\mathring{T}-\frac{1}{2}\varphi(\mathring{T})\sim (1-c)\frac{\mathring{T}}{\ln \mathring{T}}\ \Rightarrow\
\mathring{T}-T\sim (1-c)\frac{\mathring{T}}{\ln \mathring{T}}\ \Rightarrow\
\mathring{T}\sim T ,
\edis
we have (see the condition for $U$ in (1.6))
\bdis
\mathring{T}-(T+U)\sim (1-c)\frac{\mathring{T}}{\ln \mathring{T}}-U>(1-c-\epsilon)\frac{T}{\ln T}-U>(1-c-2\epsilon)\frac{T}{\ln T},
\edis
i.e. $\mathring{T}>T+U$. Then we have
\bdis
[T,T+U]\bigcap [\mathring{T},\widering{T+U}]=\emptyset;\ T+U<\mathring{T} ,
\edis
\be \label{e1.7}
d\left\{ [T,T+U];[\mathring{T},\widering{T+U}]\right\}>(1-c-2\epsilon)\frac{T}{\ln T}\to\infty ,
\ee
where $d$ denotes the distance of corresponding segments (compare \cite{11}, (1.3), (1.6)).

\begin{remark}
Some nonlocal interaction of the functions
\bdis
Z\left[\frac{\varphi(t)}{2}+\rho_1\right],\ Z\left[\frac{\varphi(t)}{2}+\rho_2\right],\ Z\left[\frac{\varphi(t)}{2}+\rho_3\right],\ \hat{Z}^2(t)
\edis
is expressed by formulae (1.4). Such the interaction is connected with two collections of disconnected sets unboundedly receding each from other
(see (1.7), $d\to\infty$ as $T\to\infty$) - like mutually receding galaxies (the Hubble law). Compare this remark with the Remark 3 in \cite{11}.
\end{remark}

\begin{remark}
If
\begin{eqnarray}
& &
x_1\not=x_2,\ x_1+x_2+2(\rho_1+\rho_2+\rho_3)\ln P\not=(2k+1)\pi, \nonumber  \\
& & \ \\
& &
y_1\not=y_2,\ y_1+y_2+2(\rho_1+\rho_2+\rho_3)\ln P\not=(2l+1)\pi, \nonumber
\end{eqnarray}
where
\be\label{e1.9}
k,l=0,\pm 1,\pm 2,\dots ,\pm L,\ L=\mcal{O}(T^{1/48-\epsilon}\ln T) ,
\ee
then the formulae (1.4) are the asymptotic formulae.
\end{remark}

\begin{remark}
The formulae (1.4) are the first asymptotic formulae (see (1.8), (1.9)) in the theory of the Riemann zeta-function for the fifth order expression
\bdis
Z\left[\frac{\varphi(t)}{2}+\rho_1\right]Z\left[\frac{\varphi(t)}{2}+\rho_2\right]
Z\left[\frac{\varphi(t)}{2}+\rho_3\right]\hat{Z}^2(t) .
\edis
These formulae cannot be obtained by methods of Balasubramanian, Heath-Brown and Ivic (see, for example, \cite{1}).
\end{remark}

This paper is a continuation of the series of paper \cite{7}-\cite{11}.

\section{Splitting of a pair of asymptotic formulae into four asymptotic formulae}

\subsection{}

In the case
\bdis
\rho_1=0,\ \rho_2=\rho_3=\rho,\ \rho=\rho_k(z)=\frac{2k\pi+z}{2\ln P},\ z\in [0,\pi]
\edis
(see (1.8), (1.9)) the formulae follows
\begin{eqnarray}
& &
\int_{\mathring{G}_5(x_1,x_2)}Z\left[\frac{\varphi(t)}{2}\right]Z^2\left[\frac{\varphi(t)}{2}+\rho_k(z)\right]\hat{Z}^2(t){\rm d}t\sim \nonumber \\
& &
\frac{4}{\pi}U\sin\frac{x_2-x_1}{2}\cos\left(\frac{x_1+x_2}{2}+z\right), \nonumber \\
& & \ \\
& &
\int_{\mathring{G}_6(y_1,y_2)}Z\left[\frac{\varphi(t)}{2}\right]Z^2\left[\frac{\varphi(t)}{2}+\rho_k(z)\right]\hat{Z}^2(t){\rm d}t\sim \nonumber \\
& &
-\frac{4}{\pi}U\sin\frac{y_2-y_1}{2}\cos\left(\frac{y_1+y_2}{2}+z\right), \nonumber
\end{eqnarray}
from (1.4). Hence, we obtain
\begin{cor}
The splitting of the formulae (2.1); $z=0,1$ leads up to the four asymptotic formulae
\begin{eqnarray*}
& &
\int_{\mathring{G}_5(x_1,x_2)}Z\left[\frac{\varphi(t)}{2}\right]Z^2\left[\frac{\varphi(t)}{2}+\rho_k(0)\right]\hat{Z}^2(t){\rm d}t\sim
\frac{4}{\pi}U\sin\frac{x_2-x_1}{2}\cos\frac{x_1+x_2}{x} , \\
& &
\int_{\mathring{G}_5(x_1,x_2)}Z\left[\frac{\varphi(t)}{2}\right]Z^2\left[\frac{\varphi(t)}{2}+\rho_k(\pi)\right]\hat{Z}^2(t){\rm d}t\sim
-\frac{4}{\pi}U\sin\frac{x_2-x_1}{2}\cos\frac{x_1+x_2}{x} , \\
& &
\int_{\mathring{G}_6(y_1,y_2)}Z\left[\frac{\varphi(t)}{2}\right]Z^2\left[\frac{\varphi(t)}{2}+\rho_k(0)\right]\hat{Z}^2(t){\rm d}t\sim
-\frac{4}{\pi}U\sin\frac{y_2-y_1}{2}\cos\frac{y_1+y_2}{x} , \\
& &
\int_{\mathring{G}_6(y_1,y_2)}Z\left[\frac{\varphi(t)}{2}\right]Z^2\left[\frac{\varphi(t)}{2}+\rho_k(\pi)\right]\hat{Z}^2(t){\rm d}t\sim
\frac{4}{\pi}U\sin\frac{y_2-y_1}{2}\cos\frac{y_1+y_2}{x} ,
\end{eqnarray*}
where $k$ fulfills the condition (1.9).
\end{cor}

\subsection{}

Next, in the case
\be \label{e2.2}
-x_1=x_2=x,\ -y_1=y_2=y,\ \mathring{G}_5(-x,x)=\mathring{G}_5(x),\ \mathring{G}_6(-y,y)=\mathring{G}_6(y)
\ee
we obtain
\begin{cor}
If $\rho_1=\rho_2=\rho_3=0$ then
\begin{eqnarray*}
& &
\int_{\mathring{G}_5(x)}Z^3\left[ \frac{\varphi(t)}{2}\right]\hat{Z}^2(t){\rm d}t\sim \frac{4}{\pi}U\sin x , \\
& &
\int_{\mathring{G}_6(y)}Z^3\left[ \frac{\varphi(t)}{2}\right]\hat{Z}^2(t){\rm d}t\sim -\frac{4}{\pi}U\sin y ,
\end{eqnarray*}
where $0<x,y<\pi/2$.
\end{cor}

\section{Law of the asymptotic equality of the areas of the positive and the negative part of the graph of the function
$Z\left[\frac{\varphi(t)}{2}+\rho_1\right]Z\left[\frac{\varphi(t)}{2}+\rho_2\right]Z\left[\frac{\varphi(t)}{2}+\rho_3\right]\hat{Z}^2(t)$,
$t\in \mathring{G}_5(x)\cup\mathring{G}_6(x)$}

Let
\begin{eqnarray*}
& &
\mathring{G}_5^+(x)=
\left\{ t:\ t\in \mathring{G}_5(x),\
Z\left[\frac{\varphi(t)}{2}+\rho_1\right]Z\left[\frac{\varphi(t)}{2}+\rho_2\right]Z\left[\frac{\varphi(t)}{2}+\rho_3\right]>0\right\}, \\
& &
\mathring{G}_5^-(x)=
\left\{ t:\ t\in \mathring{G}_5(x),\
Z\left[\frac{\varphi(t)}{2}+\rho_1\right]Z\left[\frac{\varphi(t)}{2}+\rho_2\right]Z\left[\frac{\varphi(t)}{2}+\rho_3\right]<0\right\} .
\end{eqnarray*}
Let similar definition for the sets $\mathring{G}_6^+(x),\ \mathring{G}_6^-(x)$ hold. Let
\be \label{e3.1}
\rho_1+\rho_2+\rho_3=\frac{2k\pi}{\ln P},\ k=0,\pm 1,\pm 2,\dots ,\pm L,\ L=\mcal{O}(T^{1/48-\epsilon}\ln T) .
\ee
Then we obtain from (1.4) by (3.1) the following
\begin{cor}
\begin{eqnarray}
& &
\int_{\mathring{G}_5^+(x)\cup\mathring{G}_6^+(x)}
Z\left[\frac{\varphi(t)}{2}+\rho_1\right]Z\left[\frac{\varphi(t)}{2}+\rho_2\right]Z\left[\frac{\varphi(t)}{2}+\rho_3\right]\hat{Z}^2(t){\rm d}t\sim
\nonumber \\
& &
-\int_{\mathring{G}_5^-(x)\cup\mathring{G}_6^-(x)}
Z\left[\frac{\varphi(t)}{2}+\rho_1\right]Z\left[\frac{\varphi(t)}{2}+\rho_2\right]Z\left[\frac{\varphi(t)}{2}+\rho_3\right]\hat{Z}^2(t){\rm d}t .
\end{eqnarray}
\end{cor}

Indeed, from (1.4) (see (3.1)) we have
\bdis
0<(1-\epsilon)\frac{4}{\pi}U\sin x<\int_{\mathring{G}_5(x)}\leq \int_{\mathring{G}_5^+(x)}\leq \int_{\mathring{G}_5^+(x)\cup \mathring{G}_6^+(x)},
\edis
and similarly,
\bdis
0<(1-\epsilon)\frac{4}{\pi}U\sin x<-\int_{\mathring{G}_5^-(x)\cup \mathring{G}_6^-(x)} .
\edis
Hence, adding the formulae (1.4) (see (3.1)) we get
\bdis
\int_{\mathring{G}_5^+(x)}+\int_{\mathring{G}_5^-(x)}+\int_{\mathring{G}_6^+(x)}+\int_{\mathring{G}_6^-(x)}=o(U) ,
\edis
i.e. (3.2).

\begin{remark}
The formula (3.2) represents the law of the asymptotic equality of the areas (measures) of the figures which correspond to the positive and the negative
parts of the graph of the function
\be \label{e3.3}
Z\left[\frac{\varphi(t)}{2}+\rho_1\right]Z\left[\frac{\varphi(t)}{2}+\rho_2\right]Z\left[\frac{\varphi(t)}{2}+\rho_3\right]\hat{Z}^2(t)
\ee
with respect to the disconnected sets $\mathring{G}_5^+(x)\cup\mathring{G}_6^+(x)$, $\mathring{G}_5^-(x)\cup\mathring{G}_6^-(x)$. This is one of the
laws governing the \emph{chaotic} behaviour of the positive and negative values of the function (3.3).
\end{remark}

\section{Proof of the Theorem}

First of all, from the main lemma (see \cite{1}, (2.5)) by (1.2) we have, for example,
\be \label{e4.1}
\int_{\mathring{k}_{2\nu}(x_1)}^{\mathring{k}_{2\nu}(x_2)}f\left[\frac{\varphi(t)}{2}\right]\hat{Z}^2(t){\rm d}t=
2\int_{k_{2\nu}(x_1)}^{k_{2\nu}(x_2)}f(t){\rm d}t ,
\ee
for every integrable function $f(t)$. Next, in the case
\bdis
f(t)=Z\left[\frac{\varphi(t)}{2}+\rho_1\right]Z\left[\frac{\varphi(t)}{2}+\rho_2\right]Z\left[\frac{\varphi(t)}{2}+\rho_3\right]
\edis
we have the following $\hat{Z}^2$-transformation
\begin{eqnarray}
& &
\int_{\mathring{G}_5(x_1,x_2)}
Z\left[\frac{\varphi(t)}{2}+\rho_1\right]Z\left[\frac{\varphi(t)}{2}+\rho_2\right]Z\left[\frac{\varphi(t)}{2}+\rho_3\right]\hat{Z}^2(t){\rm d}t= \nonumber \\
& &
2\int_{G_5(x_1,x_2)}Z(t+\rho_1)Z(t+\rho_2)Z(t+\rho_3){\rm d}t, \nonumber \\
& & \ \\
& &
\int_{\mathring{G}_6(y_1,y_2)}
Z\left[\frac{\varphi(t)}{2}+\rho_1\right]Z\left[\frac{\varphi(t)}{2}+\rho_2\right]Z\left[\frac{\varphi(t)}{2}+\rho_3\right]\hat{Z}^2(t){\rm d}t= \nonumber \\
& &
2\int_{G_6(y_1,y_2)}Z(t+\rho_1)Z(t+\rho_2)Z(t+\rho_3){\rm d}t, \nonumber
\end{eqnarray}
(see (1.1)-(1.3), (4.1)). Let us remind that in the paper \cite{6} we have proved following cubical correlation formulae (see \cite{6}, pp. 24, 25)
\begin{eqnarray}
& &
\int_{G_5(x_1,x_2)}Z(t+\rho_1)Z(t+\rho_2)Z(t+\rho_3){\rm d}t=\nonumber \\
& &
\frac{2}{\pi}U\sin\frac{x_2-x_1}{2}\cos\left\{\frac{x_1+x_2}{2}+(\rho_1+\rho_2+\rho_3)\ln P\right\}+\mcal{O}(T^{13/16+\epsilon}), \nonumber \\
& & \ \\
& &
\int_{G_6(y_1,y_2)}Z(t+\rho_1)Z(t+\rho_2)Z(t+\rho_3){\rm d}t=\nonumber \\
& &
-\frac{2}{\pi}U\sin\frac{y_2-y_1}{2}\cos\left\{\frac{y_1+y_2}{2}+(\rho_1+\rho_2+\rho_3)\ln P\right\}+\mcal{O}(T^{13/16+\epsilon}), \nonumber
\end{eqnarray}
($U,\rho_1,\rho_2,\rho_3$ fulfill conditions (1.6)). Now, formulae (1.4) follow from (4.2) by (4.3).

\thanks{I would like to thank Michal Demetrian for helping me with the electronic version of this work.}

\end{document}